\input amssym.tex
\input cyracc.def

\magnification =1200
\topskip = .5truecm
\vsize = 210true mm
\hsize =153 true mm
\baselineskip=15truept
\voffset=2\baselineskip
\overfullrule = 0pt
\parindent=0pt

\font\bfxx=cmbx12 at 15pt
\font\bfx=cmbx12

\font\ninesc=cmcsc9
\font\cmsc=cmcsc10

\newfam\cyrfam
\font\fivecyr=wncyr5%
\font\sevencyr=wncyr7%
\font\tencyr=wncyr10%

\def\cyrf{\fam\cyrfam\sevencyr\cyracc}

\textfont\cyrfam=\tencyr 
\scriptfont\cyrfam=\sevencyr%
\scriptscriptfont\cyrfam=\fivecyr

\font\bfn=cmbx9


\def\napa#1{\hfill\break\noindent{\bf #1}}%
%
\def\inv#1{#1^{-1}}

\def\frac#1#2{{#1\over #2}}
\def\dfrac{\displaystyle\frac}
\def\exa#1{\par\medskip\noindent{\ninesc Examples}\ (#1)\quad}%

\def\intn#1{{\cal #1}\hbox{-}\!\!\int}

\def\intdl#1#2#3#4{\int_{#1}^{#2}#3\,{\rm d}#4}

\def\intdln#1#2#3#4#5{{\cal #1}\hbox{-}\!\!\intdl{#2}{#3}{#4}{#5}}

\def\ref#1{[{\sl #1\/}]}

\newcount\notenumber

\def\noter#1{\advance\notenumber by 1
$^{\the\notenumber}$\footnote{}{$^{\the\notenumber}$\ #1}}
\outer\def\subsectionsc#1\par{\hfill\smallskip\message{
#1}\leftline{\ninesc #1}\nobreak\smallskip\noindent } %

\outer\def\proclaimsc#1#2{\par\medbreak\noindent{\ninesc #1 \enspace}{\sl
#2}\par\ifdim\lastskip<\medskipamount \removelastskip\penalty55\medskip\fi}%

%
%

\newcount\sectnumber%
\def\sectn#1 {\advance \sectnumber by 1{\par\bigskip\noindent{\bfx
\the\sectnumber\ #1}}\par}%

\newcount\subsnumber%
\def\subsn#1 {\advance \subsnumber by 1{\par\medskip\noindent{\cmsc\the\sectnumber.\the\subsnumber\ 
#1}\quad}}%

\newcount\subsanumber%
\def\asubsn#1 {\advance \subsanumber by 1{\par\medskip\noindent{\cmsc\the\sectnumber.\the\subsanumber\ 
#1}\quad}}%

\newcount\subsbnumber%
\def\bsubsn#1 {\advance \subsbnumber by 1{\par\medskip\noindent{\cmsc\the\sectnumber.\the\subsbnumber\ 
#1}\quad}}%

\newcount\subscnumber%
\def\csubsn#1 {\advance \subscnumber by 1{\par\medskip\noindent{\cmsc\the\sectnumber.\the\subscnumber\ 
#1}\quad}}%

\def\napa#1{\hfill\break\noindent{\bfn #1} }%
\pageno = 1

\rm

\centerline{\bfxx The Search for the Primitive	}

\sectn{Introduction}

In a recent interesting paper by  Gluchoff,  \ref{G}, the development of the  theory of integration	 was
related to the needs of trigonometric	series. The paper ended with the
introduction of the Lebesgue\noter{\sevenrm  Henri L\'eon Lebesgue, 1875--1941.}	  integral	at the beginning of the
century.

 More usually the case is made for the the development of the integral  being influenced by the problem
of finding anti-derivatives, or primitives.  Again  the story could end with the integral of  Lebesgue, as  this
	is the integral of  everyday mathematics. The needs of
the search for a primitive, as  Lebesgue himself called it, \ref{L}, had produced a tool of first-rate importance for
the whole of mathematical analysis. 

In neither  of these stories was the original goal attained; the  Lebesgue integral neither calculated the coefficients  of all trigonometric series,
nor gave the primitives of all derivatives. However  interest in the two problems  waned because of the overwhelming importance of the  Lebesgue
integral. Also, it was a tool of great power in the study of Fourier series, and  did find all but the most obscure primitives.

The question of completing the story for trigonometric series will be taken up in another paper,  \ref{Bu2}. Here we
will  study  how the the search for a primitive led to the  Lebesgue integral,  show how the search ended
 within a decade of the   introduction of this integral, and discuss  work that  continues to this day, more than three centuries after the basic work of
Newton\noter{\sevenrm  Sir Isaac Newton,  1642--1727.}. 

This topic is better covered in the literature than the one discussed by Gluchoff so not all
details will be given; see for example  \ref{H; L; P}. In addition the various entries in \ref{E} are useful,  and further references are
given there.

Some technicalities are placed in the Appendix, Section 8. Terms
used  which  are explained there will be written  as\  \  {\tt closed}\ \  when they first occur.
\sectn{The  Problem}

A problem that was the source of much research into integration is the following.
\proclaimsc{The Classical Primitive Problem}{ Given that function is a derivative on a closed bounded interval find the function of
which it is the derivative.}\par

This problem is well posed, as the following  simple corollary of the {\sl Mean Value Theorem of Differentiation\/} shows; see \ref{R, p.163}.
\proclaimsc{Uniqueness Theorem}{ If $F'=0$  on $[a,b]$ then $F$ is a constant function.}\par

A function $F$ such that $F'=f$ is called {\sl a primitive\/}, or {\sl an  anti-derivative\/}, of $f$, and the Uniqueness Theorem shows that   although a
derivative   can have many primitives, any two such primitives differ by a constant. 

If then  $f$ is a derivative on the interval  $[a,b]$
it has a unique primitive  $ \tilde F $ satisfying  $ \tilde F(a) = 0$; this we will call {\sl the primitive of \/}$f$, or {\sl the anti-
derivative of $f$\/}.
Equivalently if $F$ is a primitive of $f$ then the primitive of $f$ is, 
 	$$ \tilde F(x) = F(x)- F(a),\, a\le x\le b.$$

\sectn{The Newton Integral}  

Originally there was no primitive problem since all functions were
derivatives, and every function had a derivative, except at the most obvious points---like $x^{2/3}$ at the
origin---and  finding the primitive of a derivative  was solved by a theorem of    Newton.
\proclaimsc {The Fundamental Theorem of Calculus}{Given a function on
the interval $ [a,b]$ it is the derivative of the area under its graph}\par

 If the area  under the graph of $f$ and above  the interval $[a,x]$  is $A(x)$, then this theorem says that $A'= f$; since in addition $A(a) = 0$ it
follows that
$A$ is the primitive of
$f$ and the classical primitive problem is solved.

There is a minor problem:  by talking  of the area under the  graph  we seem to be assuming $f$ is positive. This can be taken care of, as in any
elementary calculus course, by understanding area to mean the area above the axis minus the area below the axis.

 Another problem was the concept of area; it  was not
defined, being  considered as an obvious idea.
\smallskip
This was the state of the subject for some time after Newton, and  is  the approach that is  common in almost  all elementary calculus courses today. The
quantity $A$ is sometimes called the {\sl Newton integral of \/}$f$.

The classical primitive problem arose because of two developments.
\smallskip
(I)  A proper definition of area was given, and the  definition  was
not related  to the definition of derivative.  

(II) The concept of function broadened, and with this came the
realization that a derivative  can be  badly behaved;  it can fail to be continuous, and  it can be unbounded,  at
many points
\smallskip
These developments meant that   now there are two classes of functions; one contains all functions  which have areas under their graphs, and one contains
all derivatives. These two classes can be distinct, and even if a function is in both of these classes we must still  see if  the Fundamental Theorem of
Calculus holds with {\it ``Given a function\dots '' \/}replaced by {\it ``Given a function that is a derivative\dots''\/}

In more common usage: we define an integral (the area), and so the class of integrable functions;  are derivatives then
integrable, and if so is the indefinite integral a  primitive? Is the derivative of the indefinite integral the original
function that was integrated?

Unfortunately the answer to the first these questions was  essentially  no. The  integral introduced  by Cauchy\noter{\sevenrm Augustin-Louis Cauchy,
1789--1857.} and Riemann\noter{\sevenrm Georg
Friedrich Bernhard Riemann, 1826--1866.}was better than the  Newton
integral, it   defined area properly and   the class of integrable functions was large. Unfortunately, although many non-derivatives
could be integrated, not  all derivatives were integrable.

During the half century following the basic work by Cauchy and Riemann, the search for an integral that would solve the classical
primitive problem was pursued. However for a solution several tools had to be in place:
\smallskip
(a)  a full understanding of the nature of sets;

(b) a good definition of measure, or area, of sets;

(c) a deeper knowledge of the properties of a derivative.
\smallskip
 These  first two  requirements were met by the beginning of the twentieth century following the work of Cantor\noter{\sevenrm Georg Ferdinand
Ludwig Cantor, 1845--1918.}, Baire\noter{\sevenrm  Ren\'e-Louis Baire, 1874--1932.}, Borel\noter{\sevenrm  Emile F\'elix-Edouard-Justin
Borel, 1871--1956.} and Lebesgue. The required properties of derivatives were obtained soon after by Denjoy\noter{\sevenrm  Arnaud Denjoy, 1884
--1974.}, who then gave the first full solution of the problem.
\sectn{The Integrals of Cauchy and Riemann}

\subsn{The Riemann Integral}  Cauchy laid the foundations  for most  of  analysis and in particular he gave the first proper definition of
area as an integral,  for  continuous functions;  and he  proved the Fundamental Theorem of Calculus for his definition. 

\proclaimsc{Theorem 1}{ If a function is a continuous derivative then  it has an area under its graph and this
area  is the  primitive of the function.}\par

 The definition of area given by Cauchy was a very methodical  form of  Archimedes\noter{\sevenrm Archimedes, 287(?)BC--212
BC.} {\sl method of exhaustion\/}, \ref{Gr, pp.29--30}. 
Look at the region under graph and divide it into  $n$  slices 
and approximate the area of each slice by  a rectangle of height the value of the function at the left-hand side of the
slice, or in fact at any point in the slice; see \ref{P, p.6 }.  If as  the width of the largest slice gets small  the total area of all $n$ rectangles has a  limit,
this says that there is an area under the graph and gives  the value of this limit as the value of the area; see  \ref{Gr,
pp.140--148; H, pp.9--11;  P, pp.3--4}.

Shortly after  this  Riemann gave the modern form of the definition. His definition  made no assumption about the nature of
the function being integrated,  the class  integrable functions, the class of functions for which an area existed,
was precisely the class for which the calculation worked.    

Importantly, Riemann gave a fundamental  necessary and sufficient condition for a function to be
integrable.

Let us give Riemann's  definition as, although it is well known,  see  \ref{K, Chap.II; R, Chap.VI}, we will need it for a later development. 

 Given a closed interval $[a,b]$ then a {\sl  partition\/} of
that interval is a finite set of points, $ \varpi=\{a_0,\ldots,a_n\}$ say, such that
$a=a_0<a_1<\cdots<a_n = b$. 

The {\sl norm of a partition\/}, $ || \varpi||$, is the length of the longest interval of the partition, that is
$|| \varpi||=\max_{1\le j\le n}|a_j-a_{j-1}|$.

A {\sl tagged partition\/} of
$[a,b]$ is a partition $\varpi$, as above, together with a {\sl tag\/} in each interval,
that is  an $ \xi_j,\, a_{j-1} \le \xi_j\le a_j,\, 1\le j\le n$. A tagged   partition will be
written $
\varpi =\{a_0,\ldots,a_n; \xi_1,\ldots, \xi_n\}$.  

If $f$ is a function defined  on $[a,b]$ and if 
$\varpi$  is a tagged  partition of $[a,b]$,  form the {\sl Riemann sum\/}
$$
S=S( f;\varpi)=\sum_{j=1}^n f( \xi_j)(a_j-a_{j-1}).\eqno(1)
$$
If  there is a real number $I$  such that given any $\epsilon>0$, there is a $\delta>0$ and
for all tagged partitions $\varpi$ of norm less than  $\delta$ 
we have that
$|S-I|<
\epsilon$  then
$f$, the {\sl  integrand\/},  is said to be {\sl  Riemann, ${\cal R}$-, integrable on \/}$[a,b]$ with {\sl integral\/} equal to $ I$; 

this written,
$$
I = \intdln{R}{a}{b}{f(x)}{x} = \lim_{n\to \infty}\sum_{j=1}^n f(
\xi_j)(a_j-a_{j-1}).\eqno(2)
$$

 The  right-hand side of (2) is an abuse of notation, in that it is not sufficient to let  $n$, the
number of points in the partition, tend to infinity; we must require that $||\varpi||$ tend to zero, when in
particular,
$n\to\infty$.

The indefinite ${\cal R}$-integral is called {\sl  an ${\cal R}$-primitive\/}, and if it is zero at the left-hand endpoint of the interval, {\sl the  ${\cal
R}$-primitive\/}. 

It is important to notice  further that  all choices of tags must be allowed, and  this has a restrictive implication.

\centerline{{\tt
For  the ${\cal R}$-integral of $f$  to exist   $f$ must be bounded.}}
\smallskip

 The basic properties of this integral are well known and are developed in all books on analysis; see  for instance \ref{R, pp.184--203}. In particular
 if $f$ is continuous, or if $f$ is monotonic, or if $f$ is continuous with a finite
number of discontinuities  then {f} is integrable; see \ref{R, pp.192--193, 196}.

The  important result of Riemann mentioned above  tells us just when a function is ${\cal R}$-integrable.
\proclaimsc{Theorem 1}{ A bounded function is Riemann integrable if and only if it is continuous {\tt almost everywhere}.}\par

It is immediate from this result  and the examples discussed in 8.4 that not all derivatives are Riemann integrable.

{\tt Derivatives
need be neither bounded  nor continuous almost everywhere}. 

However  we have the Fundamental Theorem of Calculus for this integral, due to Darboux\noter{\sevenrm Jean Gaston Darboux, 1842--1917.};
\ref{H, pp.50-51; Ho, pp.484--486; K, pp.46--50; L, p.88}
\proclaimsc{Theorem 2}{If a derivative is ${\cal R}$-integrable then the ${\cal
R}$-primitive is the primitive, that is
$(
\intn{R}_a^x f )' = f(x),\, a\le x\le b$.}\par 

\subsn{The Cauchy- Riemann Integral} Suppose that a function $f$ is continuous, or just continuous almost everywhere, but not bounded on any
interval containing a point $x$, then we will call $x$ {\sl a point of unboundedness of \/}$f$.
\exa{i} The functions $ f(x)
=1/\sqrt x$, or $g(x) = \csc 1/x$ have the origin as a point of unboundedness.  In the second case the origin is also the limit of a sequence of points of 
unboundedness; it is a {\tt limit point} of the set of points of unboundedness. In both cases  the function can be defined arbitrarily at the points of
unboundedness.

 Let $S$ be the   set of all  points of unboundedness of $f$. Then $S$ is a {\tt closed set}, and the function	 is not ${\cal R}$-integrable on any interval
that meets $S$,  but is  integrable on any interval that does not  meet
$S$. 

If $S$ is empty then the function is ${\cal R}$-integrable, and  over the next century various integrals were defined that allowed for   more and more
complicated types of non-empty  sets $S$; see \ref{P, p.10}.

 Cauchy himself had already introduced an extension of his integral method that allows us to  handle  the case where $S$ is finite..
This method is well known from elementary calculus where it goes under the name of  the {\sl  improper or infinite integral\/}; see \ref{K, pp.54--55;
R, pp.221--226}. 

Let us first assume that $S$ consists just of the point $a$, the left-hand  endpoint of the interval. Then $f$ is
 not ${\cal R}$-integrable in
$[a,b]$ but is on every $[ \alpha, b],\, a< \alpha\le b$. The {\sl Cauchy-Riemann, ${\cal  CR}$-, integral\/} of $f$ on $[a,b]$ is 
$$
I =\intn{CR}_a^b\!f=\lim_{ \alpha\to a}\intn{R}_{ \alpha}^b\!f,
$$
 provided the limit exists.

Now if $f$ is a derivative that is continuous almost everywhere on $[a,b]$, bounded on  every  subinterval $[ \alpha, b],  a< \alpha< b$, but unbounded at
a, see 8.4 Examples (i) with $\alpha= \beta = 2$ say,  then it is  ${\cal  CR}$-integrable, and its primitive is the ${\cal  CR}$-primitive.
 For suppose   $F$ is a  primitive of $f$,  by Theorem~2 
$$
\intn{R}_{ \alpha}^x\!f= F(x) - F( \alpha);
$$
 so, as $F$ is continuous,
$$
 F(x) -F(a) = \lim_{ \alpha\to 	a}\big( F(x) - F( \alpha)\big)= \lim_{ \alpha\to a}\intn{R}_{ \alpha}^x\!f= \intn{CR}_a^x\!f.
$$
The same techniques apply if  the  single point of unboundedness is $b$. 

If the single point of unboundedness is  $ c,\,  a< c< b$ then the above
definitions will give $ \intn{CR}_a^c\!f$, and $\intn{CR}_c^b\!f$, and it is then natural to define
$$
\intn{CR}_a^b\!f=\intn{CR}_a^c\!f+\intn{CR}_a^b\!f.
$$
Further, it is easy to check in these cases that   if $f$ is a derivative that is continuous almost
everywhere, with one point of unboundedness then $f$  is  ${\cal  CR}$- integrable, and  the  ${\cal  CR}$-primitive is a
primitive.

This procedure, and the result,  is readily extended to the case of $S$ being  finite; say for
simplicity that the points $S$ form a partition  $\varpi$ of $[a,b]$, $\varpi=\{s_1,s_2,\ldots, s_n \}$ then  $f$
is integrable on each $[\sigma_i, \tau_i] ,\, s_i< \sigma_i<\tau_i<s_{i+1},\, 1\le i\le n-1$ and the Cauchy-Riemann integral of $f$ is 
$$
\intn{CR}_a^b\!f=\sum_{i=1}^{n}\lim_{{ \sigma_i\to s_{\scriptscriptstyle{i-1}}}\atop{ \tau_i\to s_i}}\intn{R}_{
\sigma_i}^{\tau_i}\!f,
$$
 provided all of the limits exist.

So if $f$ is a derivative that is continuous almost everywhere, and has a
finite number of points of unboundedness then its primitive is given by the  ${\cal  CR}$-primitive; that is the Fundamental Theorem of Calculus,
Theorem 2, holds if ${\cal  CR}$ replaces ${\cal  R}$ at various points in the statement of that theorem.

Finally let us note the following alternative approach to the ${\cal  CR}$-integral. 

\proclaimsc{Theorem 3}{ Suppose $f$ is continuous almost everywhere and  has a finite set, $S$,  of  points  of unboundedness. Then $f$
is
${\cal  CR}$-integrable if and only if there is a continuous $F$ such that on any $ [\alpha, \beta]$ not containing any points of $S$,
$F( \beta) - F(
\alpha) = \intn{R}_{ \alpha}^{\beta}\!f$. Then $F$ is a ${\cal CR}$-primitive of $f$ and $\intn{CR}_a^b\!f= F(b) -F(a)$. In particular $f$ is
${\cal  CR}$-integrable if
$f$ is a derivative, and  then $F$ a primitive of $f$. 
}
 \par 
Theorem 3 gives a so-called {\sl descriptive\/} definition of the ${\cal  CR}$-integral, as opposed to the {\sl constructive\/} definition
developed  above; this approach seems to have originated with Dirichlet\noter{\sevenrm Johann Peter Gustav Lejeune-Dirichlet, 1805--1859.}; see
\ref{P, p.12}.

There is an  important distinction between the  Riemann integral and the Cauchy-Riemann integral.  Clearly $|f| $ is bounded and continuous almost
everywhere when $f$ has these properties, hence $|f| $ is ${\cal  R}$-integrable if $f $ is.  However this is not a property of  the  
Cauchy-Riemann integral.

\exa{i} Let us  consider a simple modification of  the example mentioned above, 8.4 Examples (i) with $\alpha= \beta = 2$;  $f=F'$ where
$$
F(x) = \cases{x^2\sin\dfrac{\pi}{2x^2} ,& if $0<x\le 1$,\cr
                      0,& $x=0$,\cr}
$$ 
Then $f$ consists of two parts one of which is continuous and bounded and the other is $\phi(x)=\dfrac{\pi}{x}\cos \dfrac{\pi}{2x^2}$. This graph of
$\phi$ has, as we approach the origin,  an infinite number of areas above the axis, $a_2, a_4,\ldots$ say, and an infinite number of areas below the
axis $a_1, a_3,\ldots$ say.  Then it can be checked that
$$
 \intn{R}_{1/\sqrt{2n+1}}^{1}\phi= \sum_{k=1}^n(-1)^ka_k,\;\hbox{\rm and so}\;\intn{CR}_{0}^{1}\phi=
\sum_{k=1}^{\infty}(-1)^ka_k
$$
However the series is a conditionally convergent; the  series  $\sum_{k=1}^{\infty}a_k,$ is not convergent; so $|\phi|$	   is not  ${\cal
CR}$-integrable; see \ref{K, pp.135--136}.

 This property is expressed by saying that the Riemann integral is an {\sl absolute integral\/}, while the Cauchy-Riemann integral is {\sl non-absolute
integral\/}. Clearly then, an integral that solves the primitive problem must be a non-absolute integral.
\sectn{Infinite  Sets of Unboundedness}
 \asubsn{ The Cauchy Scale of Integrals}    The above process gives what is called the {\sl Cauchy extension\/} of the Riemann integral. It
is important to note that  the process  can be applied to any  integral. In particular we could replace the basic integral by the  Lebesgue
integral. What we want to  do  is to replace the basic  integral by the  ${\cal CR}$-integral itself, and so begin an induction that will allow us to handle
a derivative with its set of points of unboundedness, $S$,  infinite but countable; that is when $S$ is {\tt reducible}.

Note that to say $S$ is finite  equivalent to saying that  to saying  that, the {\tt derived set}, $S'$, is empty. The next move would be then to
consider the case of $S$ having a finite number of limit points,  $S'$ finite, or  the {\tt second derived set}, $S^{(2)}$, empty

 Assume for simplicity that $S$ has one limit point, and  that it is $a$; see for instance the function $g $ in 3.2 Examples (i).

Then $f$ is
 not integrable on
$[a,b]$. Any  $[ \alpha, b],\, a< \alpha\le b$,   contains only a finite number of points  of $S$ and so $f$ could be  ${\cal 
CR}$-integrable on  $[ \alpha, b]$, as it would be if $f$ were a derivative. If this is the case then the {\sl Cauchy-Riemann integral  of order two \/} of
$f$ on
$[a,b]$ is 
$$
I ={\cal CR}^{(2)}\!\hbox{-}\!\!\int_a^b\!f=\lim_{ \alpha\to a}\intn{CR}_{ \alpha}^b\!f,
$$
 provided the limit exists.

By an argument analogous to that used in the case of the  ${\cal  CR}$-integral we can define the ${\cal CR}^{(2)}$- integral when $S'$ is any
finite set.

Theorem  3 can be extended to  this case:

\proclaimsc{Theorem 3$_2$}{ Suppose that $f$ is continuous almost everywhere and the set of points, $S$, of points of unboundedness has a finite
number of limit points. The function 
$f$ is
${\cal  CR}^{(2)}$-integrable if and only if there is a continuous $F$ such that  on any $ [
\alpha, \beta]$ not containing any points of  S, and  $F( \beta) - F( \alpha) = \intn{R}_{ \alpha}^{
\beta}\!f$. Then $F$ is a ${\cal CR}^{(2)}$-primitive of $f$ and ${\cal CR}^{(2)}\!\hbox{-}\!\!\int_a^b\!f= F(b) -F(a)$.
If further if $f$ is a derivative then $f$ is ${\cal  CR}^{(2)}$-integrable, and   $F$ is a primitive of $f$.}\par 

It is now an easy induction to extend this to the case when for some $n$ the {\tt  derived set of order }$n$, $ S^{(n)}$, is empty. In this way we obtain a
sequence of extensions of the Riemann integral, the     ${\cal CR}^{(n)}$-integrals, $n=1,2,\ldots$ ; where ${\cal CR}^{(1)}$ is just ${\cal
CR}$.  Further if
$f$ is a derivative, continuous almost everywhere with      $ S^{(n)}$  empty  it is    ${\cal CR}^{(n)}$-integrable and the ${\cal
CR}^{(n)}$-primitive is a primitive of $f$.

Suppose now that no   $ S^{(n)}$ is empty, then $ S^{(\omega)}=\bigcap_{n=1}^{\infty} S^{(n)}$  is also not empty, by the {\tt Cantor Intersection
Theorem}, see 8.3. 

 Assume for simplicity that
$S^{(\omega)}$  contains one point, see 8.1 Examples (v),  and that this point   is 
$a$. Then, using  the Cantor's Intersection Theorem again,  each 
$[
\alpha, b],\, a<
\alpha\le b$ fails to meet some $ S^{(n)}, n=n( \alpha)$, and so we can evaluate ${\cal CR}^{(n)}\mkern-4mu-\mkern-8mu\int_a^b\!f$,
if this integral exists and then define {\sl Cauchy -Riemann integral  of order $\omega$ \/} of
$f$ on
$[a,b]$ as
$$
I ={\cal CR}^{(\omega)}\mkern-4mu-\mkern-8mu\int_a^b\!f=\lim_{ \alpha\to a}{\cal
CR}^{(n)}\mkern-4mu-\mkern-8mu\int_{ \alpha}^b\!f,
$$
 provided the limit exists.

Proceeding with the induction we get finally a case where the derived set is empty,  see 8.1,  and the computation has been completed. 
 Call the
integral that covers all these cases the  ${\cal CR}^{(\Omega)}$-integral, and we have that  the Fundamental Theorem of Calculus, Theorem 2,
holds if ${\cal CR}^{(\Omega)}$ replaces ${\cal  R}$ at various points in the statement of that theorem.
\proclaimsc{Theorem 4}{If a derivative is  ${\cal CR}^{(\Omega)}$-integrable, that is if it continuous almost everywhere with a countable set of
points of unboundedness,  then the
${\cal CR}^{(\Omega)}$-primitive is the primitive, that is
$({\cal CR}^{(\Omega)}\mkern-4mu-\mkern-8mu\int_a^x\!f)' = f(x),\, a\le x\le b$.
}\par

However not all continuous almost everywhere derivative are ${\cal CR}^{(\Omega)}$-
integrable; see 8.4 Examples  (iii).

{\tt A derivative can be continuous almost everywhere and have an uncountable set of points of unboundedness}.

\asubsn{The Harnack Integral}
If  a derivative is continuous almost everywhere and has an uncountable set of points of unboundedness, $S$, then if $[a_n, b_n]$ is any {\tt contiguous
interval} of the {\tt perfect kernel} of $S$,  it contains only a countable subset of $S$.  Hence we can, using the ${\cal CR}^{(\Omega)}$-integral,
calculate the primitive on $[a_n, b_n]$. Since the derivative is  continuous almost everywhere the perfect kernel  will be of {\tt measure zero}, and
with the hindsight of the  Lebesgue integral we expect this set to contribute nothing to the value of the primitive. In which case we would naturally
expect the primitive to be defined by an  integral introduced by  Harnack\noter{\sevenrm A Harnack, 1851--1888.}. 
 $$
I=\intn{HR}_a^b\!f = \sum_{n=1}^{\infty}{\cal CR}^{(\Omega)}\mkern-4mu-\mkern-8mu\int_{a_n}^{b_n}\!f.\eqno(3)
 $$
This definition is valid provided the  series converge. More precisely, since the ordering of the contiguous intervals is arbitrary we will need
the series  in (3) to converge absolutely; see \ref{G, p.89; Ho vol.I, p.350; P, p.21}.

Unfortunately there are derivatives that are continuous almost everywhere  for which the series above do not converge; see \ref{D}.

In general we would expect that the {\sl Harnack-Riemann, ${\cal HR}$-, integral \/}  to be defined if (i) the ${\cal CR}^{(\Omega)}$-integral has
been computed on all the contiguous intervals $[ a_n, b_n], n=1,2,\ldots$ of a perfect set $P$,  (ii) the series of these integrals converges absolutely,
and (iii) we know what we mean by the integral over
$P$ and then
$$
I=\intn{HR}_a^b\!f =\intn{R}_P\!f+ \sum_{n=1}^{\infty}{\cal
CR}^{(\Omega)}\,\hbox{-}\!\int_{a_n}^{b_n}\!f=\intn{R}_a^b\!1_Pf+
\sum_{n=1}^{\infty}{\cal
CR}^{(\Omega)}\,\hbox{-}\!\int_{a_n}^{b_n}\!f,
$$

 where $1_P$ is the {\sl indicator function of the set $P$\/}, that is
$$
1_P(x) =
\cases
{ 1,& if $x\in P$,\cr
0,,& if $x\notin P$.\cr
}
$$
Then we have the following result due to  Lebesgue, a Fundamental Theorem of Calculus for the ${\cal HR}$-integral;	  see \ref{L, pp.209--211}.
\proclaimsc{Theorem 5}{If the derivative $f $ is ${\cal HR}$-integrable, then the  ${\cal HR}$-primitive is the
primitive, that is $\left( \intn{HR}_a^x\!f\right)'= f(x)$.}\par
\sectn{The  Lebesgue Integral}

It was natural to attempt to  give area  a direct definition and connect this with the Riemann	 integral. This was done by
Jordan\noter{\sevenrm  Marie Ennemond Camille Jordan, 1838--1922} who defined the {\sl Jordan content\/} of a set by a method that is
based on the idea of upper and lower Riemann sums; see \ref{R, pp.184--188}.

 This approach was generalized by  Lebesgue whose definition of measure for open and closed sets is given in 8.2.

 The definition of  integral that follows from this different approach is the  Lebesgue, ${\cal L}$-, integral, the basic tool of modern analysis.
The class of integrable functions  is based on a class that is naturally connected to the measure introduced, the class of {\sl measurable functions\/}. This
important vector space of functions  includes all  continuous functions, and all monotonic functions; and  importantly, it is
closed to pointwise limits, that is  the limits of sequences of measurable functions are also measurable.

The  Lebesgue integral, is well covered in all analysis texts so will not be given here; see \ref{B-B-T; Go; H; Ho vol.I; K; L; N vol.I, pp.116--184}. A few
facts are worth noting for our purpose.
\proclaimsc{Theorem 6}{(a) The   ${\cal  L }$-integral extends the  ${\cal  R }$-integral; that  is to  say it integrates every Riemann integrable
function
$f$, and then 
$\intn{ L}_a^x\! f= \intn{R}_a^x\! f,\, a\le x\le b$.\hfil\break
\phantom{aaaaaaaaaaa}(b)  A non-negative function is ${\cal  L }$-integrable if and only if the area under the graph is measurable; and the measure of
this area is equal  to the value of the integral.\hfil\break
\phantom{aaaaaaaaaaa}(c) If $f$ is measurable and  bounded  it is ${\cal  L }$-integrable.\hfil\break
\phantom{aaaaaaaaaaa}(d) The   ${\cal  L }$-integral is an absolute integral.\hfil\break
\phantom{aaaaaaaaaaa}(e) If $f=g$ almost everywhere then $f$ is ${\cal  L }$-integrable if and only if $g$ is ${\cal  L }$-integrable.
}\par
Because of  (d) the ${\cal  L }$-integral cannot solve the primitive problem, see the remark at the end of 3.2 and  Examples 3.2(i);  more generally  see
8.4 Examples (i). 

However in spite of the last remark the implications of the ${\cal  L }$-integral for our problem are significant. Derivatives are measurable, since they
are  limits of sequences of continuous functions see, 8.4, and the
following  form of the Fundamental Theorem of Calculus, due to  Lebesgue,  \ref {L p.174; P, p. 68},  is a simple consequence of the bounded convergence
theorem; see
\ref{H, p. 128}.

\proclaimsc{Theorem 7}{If a derivative is bounded then  it is  ${\cal  L }$-integrable and the ${\cal L}$-primitive is the primitive,
that is
$(\intn{L}_a^x f)' = f,\, a\le x\le b$. }\par

In one swoop the primitive problem is solved for all bounded derivatives.   In fact  Lebesgue shows that the condition of boundedness in Theorem 7, can be
replaced by ``bounded above'' or ``bounded below'', see \ref{L, p.183};  so that the only bad points that remain for our problem  are those that are
unbounded in both directions.

In addition we have the following generalization of Theorem 2, \ref{Ho vol I, pp.596--605; Ru, pp.168--169; L, p.88}
\proclaimsc{Theorem 8}{If a derivative is ${\cal L}$-integrable then the ${\cal
L}$-primitive is the primitive, that is
$(
\intn{L}_a^x f )' = f(x),\, a\le x\le b$.}\par

W can easily handle countable sets of
unboundedness by using the Cauchy extension of the  Lebesgue integral. Further the Harnack extension can be defined as for the Riemann integral leading
to  extensions of Theorems 4 and 5; see \ref{L, pp.209--211}. 
\proclaimsc{Theorem 9}{(a) If a derivative has  a countable set of
points of unboundedness,  then  it is ${\cal CL}^{(\Omega)}$-integrable, and the 
${\cal CL}^{(\Omega)}$-primitive is the primitive, that is
$({\cal CL}^{(\Omega)}\,\hbox{-}\!\int_a^x\!f)' = f(x),\, a\le x\le b$.\hfil\break
 (b) If the derivative $f $ is ${\cal L}$-integrable on the perfect kernel $P$ of the set $S$ of its points of unboundedness, and if
the  series  $\sum_{n=1}^{\infty}{\cal CL}^{(\Omega)}\mkern-4mu-\mkern-8mu\int_{a_n}^{b_n}\!f$ converges absolutely, where
$[a_n,b_n],\, n=1,\ldots$ are the contiguous intervals of $P$, that is if $f$ is  ${\cal HL}$-integrable, then the  ${\cal HL}$-primitive is the
primitive, that is $\left( \intn{HL}_a^x\!f\right)'= f(x)$.
}\par
 Unfortunately neither of the conditions in (b) need hold so again the primitive problem is not completely solved. This is almost the point that  Lebesgue 
reached in his search for the primitive. He did observe that although the conditions need not hold  on the complete interval they always do hold  on some
sub-interval that meets the perfect kernel
$P$. This crucial observation was the basis of the complete solution given by  Denjoy, see 6.1 and \ref{H, p 137}..

There is an other property of the  Lebesgue integral that on the one hand gives a descriptive definition analogous to that of the ${\cal  CR
}$-integral	 in Theorem 3, and on the other hand  provides, and solves, an interesting variant of the classical primitive problem.   

 A function $F$ is
{\sl absolutely continuous on\/}
$[a,b]$ if:
\proclaimsc{}{
for all $\epsilon>0$ there is a $ \delta>0$ such that
$$
\sum_{k=1}^n (d_k-c_k)\,<\,\delta \quad \Longrightarrow \;\sum_{k=1}^n|F(d_k) - F(c_k)|\,<\,
\epsilon,\eqno(4)
$$
where $\{[c_k, d_k],1\le k\le n\}$ is any set of non-overlapping sub-intervals of $[a,b]$;  that
is, the sum of all the $|F(d) - F(c)|$ is  to be small whenever the
sum is taken  over any collection of non-overlapping intervals $ [c,d]$ of total small length. 
}\par

 This class of functions is very important; in
particular  an absolutely continuous function has a derivative	 almost everywhere, and  we have the following extension of the uniqueness theorem.
\proclaimsc{Uniqueness Theorem$_2$}{ If  $F$ is $AC$ and $F'=0$ almost everywhere  on $[a,b]$ then $F$ is a constant function.}\par 
In addition we have the following theorem that gives a descriptive definition of the   ${\cal  L}$-integral, and solves the variant of the primitive problem
mentioned above.
\proclaimsc{Theorem 10}{(a)  A measurable function $f$ is ${\cal  L}$-integrable if and only if there is an absolutely  continuous  function $F$
such that $F'=f$ almost everywhere; and then
$\intn{L}_a^b \!f = F(b) - F(a)$. \hfil\break
\phantom{aaaaaaaaaa} (b) If $f $ is the derivative almost everywhere of an absolutely continuous function then the primitive of $f$  is the ${\cal 
L}$-primitive of
$f$. }\par

\sectn{Solutions to the Classical Primitive Problem}

\bsubsn{Denjoy's Solution}
The first solution of the classical primitive problem was given by Denjoy in 1912 by a 
method he called {\sl totalization\/}. The totalization procedure is too complicated to give in detail.  No new methods were needed  as totalization is a
combination of the  Lebesgue integral, and its  Cauchy and  Harnack extensions. What was needed however to make these simple classical
tools work was a very exhaustive study of the properties of derivatives. This was done by Denjoy, in a series of very long papers collected
together in \ref{D}. 

Totalization  depends on the following facts to get started:

(i)  for a derivative the set $S$ of points of discontinuity is a {\tt nowhere-dense} closed set; this means we can start the Cauchy extension
on the contiguous intervals of this closed set;

(ii) the set of points at which the series in the Harnack extension diverges is also  a nowhere-dense closed set; this means we can also obtain
the Harnack extension in lots of places.

As a result, after a countable set of calculations we have a contribution to the total, the integral, on all the contiguous intervals of a
nowhere-dense perfect set, $P$ say.  However   both statements (i)  and (ii) hold for every perfect set, and so we can now start again but now
using $P$ and get  the calculations on the contiguous intervals of a nowhere-dense perfect subset of $P$. After doing
this a countable number of times we arrive, using the {\tt Cantor-Baire Stationary Principle}, see 8.3,   at an empty set and the calculations are
complete.

By  earlier remarks the Denjoy integral is designed to cope with functions that are unbounded in both directions at a large set of points so that   Denjoy's
totalization procedure is equivalent to obtaining a particular sum of a non-absolutely convergent series of areas. It is a very subtle way of adding
together partial sums from the  contributions made by the various Cauchy and Harnack extensions in a particular order,  chosen to give  to give
the correct total for derivatives. The resulting integral is called the {\sl  restricted Denjoy integral\/}, the {\sl ${\cal
D}^*$-integral\/}, and every derivative is  ${\cal D}^*$-integrable, and the  ${\cal D}^*$-primitive  is a primitive.
\proclaimsc{Theorem 11}{(a)  The   ${\cal  D}^*$-integral extends the ${\cal L}$-integral, the ${\cal  CL}^{(\Omega)}$-
integral, and the  ${\cal  HL}$-integral\hfil\break
\phantom{aaaaaaaaaaa}(b) The Cauchy and Harnack extensions of  the  ${\cal  D}^*$-integral are equivalent to the  ${\cal 
D}^*$-integral.\hfil\break
\phantom{aaaaaaaaaaa}(c)  Every derivative is ${\cal  D}^*$-integrable and the primitive of a derivative is the ${\cal  D}^*$-primitive.}\par
 This last result (b)  means that in some sense we have with this integral reached the end of the road, no further generalization along the lines being
pursued are possible.
\bsubsn{ Perron's Solution}
Denjoy's solution of the classical primitive problem created a lot of interest at the time  and led
to other and simpler solutions.
The easiest of these other methods was due to Perron\noter{\sevenrm  Oskar Perron, 1880--1975.}.   

In order to
solve the problem of finding a $F$ such that  $F' = f$ on $[a,b]$ with $F(a) = 0$, consider the two
classes of functions $M, m$  defined by:
$$
 \underline DM\ge f\ge \overline Dm\qquad M(a) = m(a) = 0;
$$
(here, as elsewhere, the underbars and  overbars indicate lower and upper derivates respectively).

These classes are  called {\sl major and minor functions of\/}  $f$ respectively. The {\sl Perron integral\/} 
 is the common value of the $\inf M(b)$, and $\sup m(b)$, if such a common value exists.

It is easy to see that this integral solves the primitive problem since if $f$ is a derivative of $F$ then $F$ is both a major and a minor
function and so is the Perron integral of $f$.

This method of Perron has  important applications for  differential equations other than the present one, $y'=f$; see \ref{E, vol.7, pp.133--134, Vol.9, p
352; W}.
\bsubsn{Luzin's Solution}
Another approach, given by Luzin \noter{\sevenrm {\sevencyr
Nikola}{\sevenrm\u{\sevencyr i}} {\sevencyr Nikola}{\sevenrm\u{\sevencyr i}}{\sevencyr evich Luzin},
1883--1950; also transliterated a Lusin} within a few months of the Denjoy's announcement of his  results,  is  implicit in Denjoy's work. 

 The idea is to give a generalization of the concept of absolute continuity so that with this generalization Theorem 10 above holds for the ${\cal
D}^*$-integral.
 
The
right generalization of absolute continuity  is called {\sl restricted generalized absolute 
continuity\/}, $ ACG^*$,  which we now explain.

(a) {\it Firstly the   term `generalized, the `G' in $ ACG^*$.\/}

 Given a  property that holds on a set, then the generalized property  holds on a closed interval  if the
interval is the union of a sequence of sets on each of which the  the property holds. For example every finite function $g$  on
$[a,b]$ is generalized bounded since $[a,b]$  is the union of the sets on which $|g|\le n$. This simple idea is behind the extension; the 
Lebesgue inverts bounded derivatives, the ${\cal D}^*$-integral  inverts all finite derivatives, that is generalized bounded derivatives; see \ref{ Bu1,
p.458}.

(b) {\it Now to extend the definition of  absolute continuity on an interval to 
 absolute continuity on a set\/}.

To define {\sl $AC$ on $E$\/},   all that has to be done   is to
require that the endpoints in (4), $\{c_k, d_k, 1\le k\le n\}$, lie in the set $E$.

 Then we get the class  $ACG$ by letting the
interval be the union of a sequence of sets on each of which the function is $AC$.

However  this class is too  general, as a continuous ACG function need not have a derivative almost everywhere; \ref{Go, p.101; S,
p.224}.

(c) {\it Finally the meaning of the term  restricted, or  the  star in   $ACG^* $ }.

There is another form of  definition (4) in  which the quantities 
$|F(d_k) - F(c_k)|$   are replaced by
$\omega	(F;[c_k,d_k])$, the {\tt oscillation of the function $F$ on the interval $[c_k,d_k]$}. These two definitions are equivalent on
intervals, see \ref{Ho vol.I, pp. 331--337}, but   a moment's reflection will show that this is not the case for sets, as the use of 
$\omega$ involves values of $F$ off the set $E$. If the oscillation form   is used, we say that a function is absolutely continuous in the narrow or
restricted sense,  is $AC^*$, on a set; this is a stronger requirement than AC.  Then we get the  class $ ACG^*$ by letting the interval be he union
of a sequence of sets on each of which the function is
$ AC^*$.

\exa{i} The function
 $$
1_{\Bbb Q}(x) = \cases{1, & if $x\in \Bbb Q\cap[0,1]$,\cr
                              0, & if $ x\in [0,1]\setminus\Bbb Q$,\cr
}
$$
 is AC on both sets $\Bbb Q\cap[0,1],  [0,1]\setminus\Bbb Q$,
being constant on each,  and so  is ACG on $[0,1]$; but is not $AC^*$ on
either set, and is not $ ACG^*$ on $[0,1]$.

$ACG^*$functions have derivatives almost everywhere and we have the following extension of the uniqueness theorem.
\proclaimsc{Uniqueness Theorem$_3$}{If  $F$ is continuous, $ACG^*$ and $F'=0$ almost everywhere  on $[a,b]$ then $F$ is a constant
function.}\par 
In addition we have the following theorem that generalizes Theorem 10,  gives a descriptive definition of the   ${\cal  D}^*$-integral, and gives another 
variant of the primitive problem.
\proclaimsc{Theorem 12}{(a)  A measurable function $f$ is ${\cal  D}^*$-integrable if and only if there is a continuous $ ACG^*$ function $F$
such that $F'=f$ almost everywhere; and then
${\cal D}^*\mkern-7mu-\mkern-8mu\int_a^bf = F(b) -
F(a)$. \hfil\break
\phantom{aaaaaaaaaaa} (b) If $f $is the derivative almost everywhere of a continuous $ACG^*$ function then the primitive of $f$  is the ${\cal
D}^*$-primitive of
$f$. }\par

A final remark about the $ ACG$ class of functions, which are in some sense the widest general class of interest in analysis. Such functions need
not have derivative almost everywhere, see \ref{Go, p.101; S, p.224};  but  if an an $ ACG$ does have a derivative almost everywhere then we get yet
another variant of the primitive problem.
\proclaimsc{Uniqueness Theorem$_4$}{ If  $F$ is continuous and $ACG$, and has a derivative almost everywhere  with  $F'=0$ almost
everywhere  on
$[a,b]$ then
$F$ is a constant function. }\par 
A function $f$ is then said to be Denjoy-Hin\v cin\noter{\sevenrm{\cyrf  A   Ya
Hinchin}, 1894--1959; also transliterated as  Khintchine.}, or ${\cal DH}$-,(or sometimes 
${\cal DK}$-, see footnote 16), integrable, with
$F$ a
${\cal DH}$-primitive if  there is an $F$,  continuous, differentiable almost everywhere, $ACG$ and with $F'=f$ almost everywhere; see
\ref{Go; S}. .
\bsubsn{The Generalized Riemann Integral}
Another approach to the primitive problem was introduced much later by  Kurzweil\noter{\sevenrm  Jaroslav Kurzweil,1926--.} and
Henstock\noter{\sevenrm Ralph Henstock, --.} independently; see [Go, Chap.9; He; Ku]. 
Their method   generalizes the 
Riemann approach to integration by replacing the uniformly small partitions of intervals used 
in the elementary definition, see 3.1,  by partitions that are
locally small.  

Given a positive function $\delta$, a tagged partition $\varpi=\{a_0,\ldots,a_n; \xi_1,\ldots, \xi_n\}$ is said to be {\sl $\delta$-fine\/}
if
$$
 \xi_i-\delta( \xi_i)<x_{i-1}<  \xi_i<x_i<  \xi_i
+\delta ( \xi_i),\; 1\le i\le n.
\eqno(\delta)
$$
Then  we say that
$I$ is the {\sl Henstock-Kurzweil integral, $\cal HK$-integral\/},	  of
$f$  on
$[a,b]$ if given any $ \epsilon>0$ there is a positive function $\delta(x)$ such that for all $ \delta$-fine
partitions of the interval $[a,b]$  we have that $|S-I|< \epsilon$, $S$ as in (1).
\exa{i} If $f$ is the derivative of $F$, and $\epsilon>0$ define the positive function $ \delta$ by
$$ |F(v)-F(u)-hF'(x)|=
|F(v)-F(u)-(v-u)f(x)| <|v-u| \frac{\epsilon}{b-a},
$$
if
$|v-u|< \delta$, and let $\varpi$ be a $\delta$-fine partition, as above. Then, 
$$
\eqalign
{
\left|F(b) - F(a) -\sum_{i=1}^nf( \xi_i)(a_i-a_{i-1})\right|=&\left|\sum_{i=1}^n\bigl(F(a_1) - F(a_{i-1}) -(a_i-a_{i-1})f(
\xi_i)\bigr)\right|\cr
\le&\sum_{i=1}^n\left|F(a_1) - F(a_{i-1}) -(a_i-a_{i-1})f(
\xi_i)\right|\cr
\le&\frac{ \epsilon}{b-a}\sum_{i=1}^n(a_i-a_{i-1})=		\epsilon.	\cr
}
$$
So $f$ is  $\cal HK$-integrable and its primitive is a  $\cal HK$-primitive.

Simple approaches to this integral can be found in \ref{D-S; L-V; McL;McS}.

\bsubsn{ Relationships Between the Various Solutions}
 The various definition of the integrals given above are very  different and a natural question  is --- how are they related? They all extend the
${\cal L }$-integral, and in fact extend the ${\cal HL}$-integral. That is to say if a function is integrable in the  Lebesgue sense
then it is integrable	 in all the above senses and the integrals are equal; more, if the functions are non-negative or just bounded below,  then if
they are integrable in any of these senses they are also  Lebesgue integrable. In addition these integrals have no Cauchy or Harnack extension---
that is any such  extension just give back the original integral; see Theorem 11(b). Finally all solve the primitive problem, that is Theorem 11 (c) holds
with each of these integrals replacing the ${\cal D}^*$-integral

 Surprisingly all these integrals are equivalent.
\proclaimsc{The Hake-Looman-Aleksandrov Theorem}{If a function is integrable by the method of Luzin, or  ${\cal  P}$-, ${\cal  D}^*$- or
${\cal  HK}$-integrable then it is integrable by all these methods, and the integrals are equal.}\par
 
By far the most difficult part of this  result is the equivalence of the Perron and restricted Denjoy integral, a problem solved by the three
mathematicians H  Hake, H   Looman and
Aleksandrov\noter{\sevenrm{\cyrf  P S Aleksandrov}, 1896--1982; also transliterated as Alexandroff.}  about ten years after Denjoy's
initial paper.

 The theorem says that the various integrals defined  to solve the classical primitive problem are just different
ways of looking at the same integral. This need not have been the case; as we  see elsewhere the various integrals involved in the  solution to the classical
coefficient problem are not equivalent; see \ref{Bu2}.

 Details    of these integrals can be found in  the classical book of Saks, \ref{S.
pp.186--259}, as well as in the  recent book by Gordon,\ref{Go}; see also \ref{Br; N vol.II; P}. 

\bsubsn{The Power of the Solution}
 The classical primitive problem is well posed because,  if $F'$ is zero everywhere then  $F$ is a constant. There are many
generalizations of this result and we can ask whether the methods discussed will solve these generalizations. 

Perhaps the most elementary result is the following, see \ref{B, pp.19--21}. 
\proclaimsc{Uniqueness Theorem$_5$}{If $F$ is continuous and  $F_+' =0$ {\tt nearly everywhere } then $F$ is a constant function.}\par

So if $f$ is nearly everywhere the right derivative of a continuous function we can ask  for  its primitive just as in the classical primitive problem. It can
be shown that $f$ is ${\cal D^*}$-integrable, and that the primitive is the ${\cal D^*}$-primitive.

 Clearly the same problem, with same solution exists
for the left derivative.

In this problem the requirement of continuity is essential, it cannot even be replaced by right continuity\noter{\sevenrm The elaborate
counter-example given  in [B, p.27] seems to be false; see [Bus]. }.
\exa{i} Define
$$
F(x) = 
\cases{
0,& if $0\le x< 1/2$,\cr
\alpha,& if $x= 1/2$,\cr
\beta,& if $1/2<x\le1$,\cr
}
$$

This example shows that we cannot allow one exceptional point in the classical primitive problem, since for all choices of $\alpha$ and  $\beta$,  $F'(x)=
0, x\ne 1/2$. In addition if $ \alpha= \beta$ $F$ is right continuous, so justifies the above remark.

In addition we cannot allow an uncountable set exceptional set, even if it is of measure zero,  as the following example shows; see \ref{Go. pp.13--14;
Ru, p.168}.

\exa{ii} Let  $K$ be the {\tt  Cantor ternary set}.  Now define $F$ as follows: on the first removed open interval, of length $1/3$, $F$ is put equal to
$1/2$; on the two removed intervals of length $1/9$, put $F$ equal to $1/4$ on the left interval and $3/4$ on the right interval; on the four removed
intervals of length $1/27$ put $F$ equal to $1/8, 3/8, 5/8,7/8$, as the intervals are chosen from left to right; etc. In this way we define  $F$ off  the
set $K$; on
$K$ it is defined by continuity. The resulting function $F$,  called the {\sl  Cantor function\/},  is increasing, with $F(0)= 0, F(1)= 1$ and off  $K$
clearly 
$F'= 0$; that is
$F'=0$ almost everywhere.

This function	 is not of course  constant ;  in fact looking at Theorems 10 and 12 we deduce that $F$ is not $AC $, or even $ACG^*$. In fact  from the
comments at the end of 6.3 it is not $ACG $ either.

To obtain a correct extension of the primitive problem if uncountable exceptional sets are allowed we must either, as in Theorems 10 and 12 restrict the
class of primitives even more than above, to $ AC$ or $ACG^*$ functions, or have some other knowledge as in the following  theorem; see \ref{Br,
p.120; K-K, pp.102--103; S, pp.205--207}.

\proclaimsc{Uniqueness Theorem$_6$}{If $F$ is continuous and  differentiable nearly everywhere with $F'=0$ almost everywhere then $F$ is a
constant function.}\par

So if $f$ is almost  everywhere the  derivative of a continuous nearly everywhere differentiable function we can ask  for  its primitive just as in the
classical primitive problem. It can be shown that $f$ is ${\cal D^*}$-integrable, and that the primitive is the ${\cal D^*}$-primitive; see\ref{Go,
p.108}.

There are many other variants of the uniqueness theorem to be found in the various references. The discussion above shows the power of the integrals
developed in that they are able to handle most of the problems that arise. The simplest way is to show that  under the conditions given by a uniqueness
theorem the primitive will be  $ACG^*$, when the  ${\cal D^*}$-integral will give the solution. 

A rather surprising result is that recently the  ${\cal D^*}$-integral has been shown to be too powerful. It is natural to ask if there is a minimal integral
that will generalize the  ${\cal L}$-integral and integrate all derivatives; that is this integral will do just that and no less general integral will.
Surprisingly there is such an integral and it is not the  ${\cal D^*}$-integral. its construction is based on an ingenious extension of the ${\cal
H-K}$-integral.

If  ($ \delta$) is replaced by
$$
 \xi_i-\delta( \xi_i)<x_{i-1}< x_i<  \xi_i
+\delta ( \xi_i),\; 1\le i\le n.
\eqno(\delta^*);
$$
that is if we do not require the tag to lie in the sub-interval that it tags  but jut close by,  then the integral defined, often called the {\it McShane
integral\/},  is just the   ${\cal L}$-integral; see \ref{L-V p,127], McS}. If   in the
definition of the integral the distance of the tag from the subinterval  is restricted by requiring  that   $\sum_{i=1}^n {\rm dist}(\xi_i,[x_{i-1}, x_i])
<1/\epsilon$ we get an integral, called the   ${\cal C}$-integral, that is strictly more general that the  ${\cal L}$-integral  and strictly less general than
the  ${\cal D^*}$-integral  and is the minimal integral ; see \ref{B-D-P} .
				 
\sectn{ Appendix}

             Everything we have to  say will occur in a closed bounded, compact, 
interval that may or may not be specified as $[a,b]$.

The {\sl oscillation of a function $f$ on the interval $[a,b]$\/} is
 $$
\omega(f;[a,b]) = \sup_{a\le x,y\le b}|f(y)-f(x)| = \sup_{a\le x\le b}f(x)- \inf_{a\le x\le b}f(x).
$$
A property that holds except on a countable set is said to hold {\sl nearly everywhere\/}.
\csubsn{ Open and Closed Sets}
An {\sl open set\/}, $G$ say, is  either the empty set, or it is a union of a countable collection of non-overlapping  open
intervals,
$I_n=]a_n,b_n[, n=1,2,\ldots$ say; so $G= \bigcup_{n=1}^{\infty}I_n$
\smallskip
A {\sl closed set\/} is either the empty set, or it is the complement, in the basic bounded closed interval,  of an open set; so if $G$ is as above
$F=[a,b]\setminus G = [a,b]\setminus \bigcup_{n=1}^{\infty}I_n$ is a closed set.  The  closed intervals $[a_n,b_n], n=1,2,\ldots$, are
called the {\sl contiguous intervals\/} of $F$.

A point is said to be a {\sl limit point of a set\/} if there is a sequence of distinct points of the set that converge to  that point.

 A closed set is
distinguished by the fact that every limit point of the set is  in the set.  Not every point  of a closed set need be be a limit point--- if that is
the case then the set is called a  {\sl perfect set\/}. 

A finite set has no limit points. An infinite set, in $[a,b]$, must have at least one  limit point by the Bolzano-Weierstrass
Theorem\noter{\sevenrm Bernard  Placidus Johann Nepomuk Bolzano, 1781--1848.}$^,\!\!$\ \noter{\sevenrm
Karl Wilhem Theodor Weierstra\ss, 1815--1897.}; see 8.3.

\exa{i} If  $ [a,b] = [0,1]$ and $S= \{1, 1/2, 1.3,\ldots, 1/n, \ldots, 0\}$ then the only limit point is $0$; its contiguous intervals are
$[1/(n+1),1/n],\,  n=1,2,\ldots$

(ii) If $S= [a,b]$ then $S$ is closed and  every point is a limit point, so $S$ is perfect.
\smallskip
Starting with a non-empty closed set  $S= S^{(0)}$ say, let the set of its limit points be $S'=S^{(1)}$; this set is  called the {\sl first derived
set of $ S$\/} It is also a closed set,  see \ref{R, p. 55; N  vol.I, p. 37}, and is  of course a subset of $S$. Further it is a proper subset unless $S$ is
perfect.

Unless
$S'$ is empty, or equivalently $S$ is finite, this  process can be repeated  by taking all the limit points of 
$S^{(1)}$, getting the {\sl second derived set of $ S$\/}, $S^{(2)}$, a subset of $S^{(1)}$. Thus we can define  a sequence
$S^{(0)}, S^{(1)}, S^{(2)},\ldots$ of derived sets of higher and higher order.
\smallskip 

There are three possibilities: 

(i) at some stage the {\sl derived set of order \/}$n$,   $S^{(n)}$, is empty, equivalently  $S^{(n-1)}$ is finite; when of course the set $S$ is countable;

(ii)  at some stage  $S^{(n)}$ is not empty but is not  a proper subset of $S^{(n-1)}$; that is $S^{(n-1)}= S^{(n)}\ne\emptyset$. So  $S^{(n-1)}$ 
is a perfect set, all its points are limit points, and so $S$ is not countable;

(iii) neither of these happens, when not only do we get a different $S^{(n)}$  for all $n$, but in addition $S^{(\omega)}
=\bigcap_{n=0}^{\infty}S^{(n)}$  is not empty, by the Cantor Intersection Theorem, see 8.3; this can happen  both when  $S$ is countable, and 
when $S$ is uncountable;

\exa{iii} If S is as in  Example (i) $S^{(1)}=\{0\}$, $S^{(2)}=\emptyset$.

(iv) If $S = [1,2]\cup S_1$, where $S_1$ is the set in Example (i), then $S^{(1)}=\{0\}\cup [1,2]$, and $S^{(2)}=S^{(3)}=\cdots
=[1,2]$.

(v) Let $a_n, n=1,2,..$ be any strictly decreasing sequence in $[0,1]$ with limit $0$ and put $A_1=\{0,
a_1,a_2,\ldots\}$. Now on  left-hand half of  each $[a_{n+1}, a_n]$ put a copy of $A_1$, and call the new set $A_2$.
Repeat this on each interval between two elements of $A_2$. Keep on doing this and call the final set, the union
of all this procedure,
$S$. Then  case (iii) above occurs, $ S^{(n)}$ exist, and is non-empty,  for all $n$, and is a proper subset of $S^{(n-1)}$; $S^{(\omega)}
=\{0\}$.  Clearly $S$ is countable
but $S\cup[1,2]$ is uncountable, and this illustrates the last remark in case (iii).
\smallskip
Suppose then $S$ is such that case (iii) holds. Then  we can start the process all over again  starting with $ \tilde S =S^{(\omega)}$ to get $ \tilde
S'=S^{(
\omega+1)}, S^{(
\omega+2)}\ldots$. This  new sequence  can exhibit all  the three possibilities  mentioned above. If again case (iii) occurs then the process can be
repeated starting with  $S^{( 2\omega)}=\cap_{n=1}^{\infty}S^{( \omega + n)}$.
\exa{vi} If $S$ is Examples (v)  $S^{( \omega+1)}$ is empty.

(vii) If $S$ is Examples (v) to together with  $[1,2]$, then $S^{ (\omega+1)}=\cdots =[1,2]$.
\smallskip

 It is a very important result   that sooner or later either  case (i) or (ii) must arise; see \ref{Ho vol.I, pp.124--125}. This is an example of the
important  Cantor-Baire Stationary Principle; see 8.3.

 If the above procedure terminates with case (i)  the original set was countable, and is called
a {\sl reducible\/} set. Otherwise the  the set was uncountable  and is  said to be {\sl irreducible\/}; the resulting perfect set is called the {\sl perfect
kernel\/}, or {\sl nucleus\/}, of the original closed set.
\csubsn{Measure}
The {\sl length of any interval \/} $I$ with endpoints $a$ and $b$ is $|I|=
|b-a|$.

A set $E$ is  said to be of {\sl measure zero\/} if: 
 give any $\epsilon>0$ there are
 intervals  $I_n, n=1,2,\ldots $such that 
(i)  $E\subseteq \cup_{n=1}^{\infty}I_n$, and (ii)
$\sum_{n=1}^{\infty}|I_  n|<\epsilon$.

The complement, in our basic closed bounded interval,  of a set of measure zero will be said to be of {\sl full measure\/}.

A property that holds on a set of full measure is said to hold {\sl almost everywhere\/}.

\exa{i} Any countable set $\{c_1,c_2, \ldots\}$ has measure zero as can be seen by putting $c_n$ inside
$I_n=]c_n - \epsilon2^{-n-1}, c_n + \epsilon2^{-n-1} [, n=1,2.\ldots$. In particular the empty set,  any finite  set, and the set of rationals $\Bbb
Q$ is of measure zero.

(ii) The set of irrationals is is of full measure.
\smallskip
 The {\sl measure of  the open set\/} $G= \bigcup_{n=1}^{\infty}I_n$ is   $|G|= \sum_{n=1}^{\infty}|I_n|$ 

  The {\sl  measure of the  closed set \/}$F=[a,b]\setminus G$, is $|F|=b-a-|G|$.
\exa{ii} If the open set $G$ is of full measure then the closed set $F$ is of measure zero. This is easily checked
since for every $m$, $F$ is contained in  the finite collection of closed intervals that make up the set
$[a,b]\setminus  \bigcup_{n=1}^{m}I_n$.
\smallskip

It is important to be  be able to exhibit a closed set of every measure from zero to full; see \ref{K. pp.84--85}.
The  extremes are easy:
finite sets are of zero measure, and of course $[a,b]$ is itself of full measure These examples are not very
useful and in the case of zero measure we can do better. There are many ways of constructing closed sets of a given measure but we will do it in a
systematic and simple way.

Given $I_{01} = [a,b]$ put $\mu _0 = |I| = b-a$ and suppose given a sequence $\underline\varepsilon$ of positive real
numbers $\varepsilon_n$, with $0<\varepsilon_n<1,\, n=1,2,\ldots$. 

Remove from $I_0$ the central open interval of
length $\varepsilon_1\mu_0$, leaving two symmetrical closed intervals  $I_{11},\, I_{12}$ both of the same length,
$\mu_1$ say; and clearly  the total length of the two closed intervals is  $2\mu_1 =\mu_0 -
\varepsilon_1\mu_0=\mu_0(1-
\varepsilon_1) = (b-a)(1-
\varepsilon_1) $.

 Remove from each of these closed intervals,  $I_{11}$ and $ I_{12}$, the central open interval of length
$\varepsilon_2\mu_1$, leaving four symmetrically situated closed intervals $I_{2i},\, 1\le i\le 4$, each of
length
$\mu_2$; clearly  $2\mu_2 =\mu_1 - \varepsilon_2\mu_1=\mu_1(1- \varepsilon_2)$ so the total length of the four closed intervals is  $4\mu_2
=2\mu_1(1- \varepsilon_2)=\mu_0(1- \varepsilon_1)(1- \varepsilon_2)= (b-a)(1- \varepsilon_1)(1-
\varepsilon_2)$. 

At the $n$th stage we remove $2^{n-1}$ central open intervals  each of length
$\varepsilon_n\mu_{n-1}$ leaving $2^n$ intervals $I_{ni},\, 1\le i\le 2^n$ of length $\mu_n$.
The total length of the closed intervals remaining at the $n$th stage is
$$
2^n\mu_n = \mu_0\prod_{k=1}^n (1-\varepsilon_k)=(b-a)\prod_{k=1}^n (1-\varepsilon_k), 
$$ 
 and so the length of the removed open intervals  is  $(b-a)\left(1-\prod_{k=1}^n (1-\varepsilon_k)\right)$.

 Let then $G$ be the open set that consists of  the complete sequence of  open intervals removed in this
manner; and put $K $ equal to the closed set that is the intersection of all the remaining closed intervals obtained at each stage. Then
$K= [a,b]\setminus G$   is not empty, by the Cantor Intersection Principle, see 8.3,  and in any case is easily seen to contain at least the end points of all
the removed open intervals.  This set $K $ is often called a {\sl generalised 
Cantor set\/}, or a {\sl Cantor-like set\/}; see \ref{B-B-T, p.28}.
 
From the above construction we have that:
$$
|G| = (b-a)\left(1-\prod_{k=1}^{\infty} (1-\varepsilon_k)\right),\quad  |K| = (b-a)\prod_{k=1}^{\infty}
(1-\varepsilon_k).
$$
Now from  the elementary theory of infinite products, see  for instance \ref{Kn, pp.218--229}, $\prod_{k=1}^{\infty}
(1-\varepsilon_k)$ converges if and only if $\sum_{k=1}^{\infty}\varepsilon_k$ converges, while  the product diverges to $0$ if 
$\sum_{k=1}^{\infty}\varepsilon_k$ diverges.

Hence: if $\sum_{k=1}^{\infty}\varepsilon_k<\infty$  then $0<|K|<b-a$;
while if  $\sum_{k=1}^{\infty}\varepsilon_k=\infty$    then $|K|=0$.

\exa{i} If $\varepsilon_n = 1/3,\, n= 1,2,\ldots  , a=0, b=1$,  then 
obviously $\sum_{k=1}^{\infty}\varepsilon_k=\infty$; the above set, which is of zero measure,  is known as the
{\sl Cantor ternary set\/}, or often just {\sl the Cantor set\/}.

(ii) If $0<\theta<\pi,\,\varepsilon_n = \theta^2/(n^2\pi^2),\, n= 1,2,\ldots $ then
$\sum_{k=1}^{\infty}\varepsilon_k<\infty$, see \ref{Kn, p.221}, and so $K$ has positive measure.  In this case
the value of the infinite product is known to be $\sin\theta/\theta$, \ref{Ru, p.310}. Since this last function takes
every value between $0$ and $1$  we see that  by the right choice of $\theta$, $K$  can have every measure greater than zero and less than $(b-a)$
\smallskip
It is worth noting some other properties of the Cantor sets: 

(a) $K$  is a perfect set;

(b) $K$ is uncountable;

(c)  in every neighbourhood of a point of $K$ can be found an interval of $G$; this is expressed by saying that  $K$
is {\sl nowhere-dense\/};

(d) when $K$ has positive measure so does $G$ and further they have a kind of fractal property; in every
neighbourhood of every point of $K$ there is a part of $K$ of positive measure, and the same is true for $G$; they
are said to be {\sl thick-in-themselves\/}.
\csubsn{ Some Basic Theorems} The following results are well known but will be referred to by name and are listed here, in the forms needed,  for
convenience.
\proclaimsc{Bolzano-Weierstrass Theorem \ref{R, p.53; N vol.I, pp.35--36}.}{Every bounded closed set contains at least one limit point.}\par
\proclaimsc{Cantor's Intersection Theorem \ref{B-B-T, pp.8--10; R, p.64}.}{ If  the intersection of a decreasing sequence of bounded  closed sets is
empty then one of the sets in the sequence is itself empty}\par

 In particular this shows that if none of the sets of the decreasing sequence is empty then their intersection cannot be empty.
\proclaimsc{ Cantor-Baire Stationary Principle \ref{Br, p.55; N vol.II, p.145}.}{In the construction of a decreasing  family of bounded closed sets one
must after a countable number of steps arrive at a point where all the sets in the construction are the same.}\par
 In particular if the family is strictly decreasing this means that after a countable number of steps the members of the family must all be empty.
\csubsn{ Discontinuous Derivatives} A derivative need not be continuous, and
it need not be bounded as the following standard example shows,
\exa{i} If  $\alpha >0$ and $\beta>0$ define, 
$$
\phi_{0}(x) = \cases{x^{\alpha}\sin x^{-\beta},& if $0<x\le 1$,\cr
                      0,& $x=0$,\cr}
$$
 
 Then $\phi_0$  is continuous, and  if $\alpha > 1$\   is differentiable with 
$$
	\phi_0'(x) = \cases{\alpha x^{\alpha-1}\sin  x^{-\beta} - \beta x^{\alpha-\beta-1}\cos x^{-\beta}, & if   $0<
x\le1,$\cr
                        0,& $x= 0$.\cr
}
$$
 So $\phi_0'$
 is continuous at the origin if   $\alpha>\beta+1$, and is not continuous at the origin if $1<\alpha\le\beta+1$; it is  unbounded there if  
$1<\alpha<\beta+1$. In addition if  $1<\alpha\le\beta$ then $\phi_0'$ is not ${\cal L}$- integrable in any interval that contains the origin; see
\ref{Br, p.52; Bu1}; the standard example has $\alpha= \beta=2$; see \ref{K, pp.135--136}.
\smallskip
 
This example can be elaborated by a standard process to produce a continuous function
whose derivative has a dense set of discontinuities.

\exa{ii}  Simple modifications of $\phi_0$ will give a function defined on $[a,b]$ and which has the same characteristics  at both $a$ and $b$ as
$\phi_0$ has at the origin:
$$
\phi_{a,b}(x)=
\cases{
 \left(\dfrac{(x-a)(b-x)}{b-a}\right)^{\alpha}\sin \left(\dfrac{(x-a)(b-x)}{b-a}\right)^{-\beta},  & if $ a<x<b$,\cr
\cr
0, & if $x=a$, or  $x=b$.
}
$$
It is worth noting that $ |\phi_{a,b}|\le \left((b-a)/4\right)^{\alpha}$

(iii) If  $\phi_0$ is as in Examples (i) and $r\in\Bbb R$ write $\phi_r(x) = \phi_0(x-r)$  with $\alpha, \beta$
depending on $r$. Let $C =  \{c_1,c_2,\ldots \}$, be any countable set,  possibly dense such as the rationals, and put
$$
\phi(x) = \sum_{n=1}^\infty\frac{\phi_{c_n}(x)}{n!},\quad  x\in\Bbb R.
$$
 
Then we easily see that if $ \alpha_{c_n}>1, n= 1,2,\ldots$, $\phi$ is continuous and differentiable with
$$
\phi'(x) = \sum_{n=1}^\infty \frac{\phi_{c_n}'(x)}{n!}.\quad x\in\Bbb R;
$$
 further $\phi'$ is continuous except at the points of $C$ where, at $c_n,\,\phi'$ has the same
discontinuity as $\phi_{c_n}'$.
\smallskip

 In this way the set of points of discontinuity of a derivative can be  countable, and can also  be dense as we can take  as the
countable set to be the set of rationals; \ref{Br, p. 34}.

 For further developments it is necessary to give a more elaborate
example using the generalised Cantor set of 2.1.

\exa{iii}
Now if $K$ is a generalized Cantor  set, as in 8.2,  with  contiguous intervals
$I_n = [a_n,b_n],\, n\in\Bbb R^*$ define , using Examples (ii) above,  $\phi_K:[a,b]\to\Bbb R$ as follows
$$
\phi_K(x) =\cases{\phi_{a_n, b_n}(x),&\quad $x\in I_n,\, n=1,2,\ldots$;\cr
 0,&\qquad  $ x\in K$.\cr 
}
$$
Here the real numbers $\alpha,\, \beta$ are subject to the same conditions as  Examples (i), (ii); in particular we always assume that $ \alpha,\, \beta$
are positive and so $\phi_K$ is always continuous.  It follows that $\phi_K$ exhibits at each point of $K$ the
character that $\phi_0$ exhibited at the origin, namely:\hfil\break 
\phantom{a} if  $\alpha>1$  then $ \phi_K$  is differentiable ; \hfil\break
\phantom{a} if $\alpha >\beta +1$ then $ \phi_K '$ is  continuous ;\hfil\break
\phantom{a}  if $1<\alpha\le \beta +1$  then 
$\phi_K'$  is not continuous at any   point  in $K$;\hfil\break
\phantom{a} if $1<\alpha<\beta+1$ then every point of $K$ is a point of unboundedness of  $\phi_K'$;\hfil\break
\phantom{a} if  $ 1<\alpha\le \beta$ then $ \phi_K'$ is  not ${\cal L}$-integrable  in any
neighbourhood of  any  point  of  $ K$. 
\smallskip
This example  is due to  Volterra\noter{\sevenrm Vito Volterra,1860--1940.}; \ref{Ho, vol.1, pp.490--491; J, pp.148--149}. It is
to be noted from 8.2 that  $K$ can be so chosen that both it and its complement, $G$,  are quite thick; both sets being
in thick-in-themselves.

\exa{iv} Now if $n=2,3,\ldots $ let $K_n$ denote a generalized Cantor set with measure
greater than
$\bigl(1-\dfrac{1}{n}\bigr)(b-a)$; further let $\phi_n$ be the $\phi_{K_n}$	 constructed as in Examples (iii), the
numbers $\alpha,\, \beta$ being independent of $n$. Now define $\phi:[a,b]\to\Bbb R$ by
$$
\phi(x) = \sum_{n=1}^\infty \frac{\phi_n(x)}{3^n},\quad a\le x\le b.
$$
Then $\phi$ is differentiable and its derivative has $E = \bigcup_{n=1}^\infty K_n$ as its points of discontinuity,
unboundedness, or non ${\cal L}$-integrability, depending on the choices of $\alpha$ and $\beta$. Since now $|E| = b-a$,
we have a derivative with a set of points of discontinuities, etc., that has full
measure. A further discussion of this  can be found in the references.
 
Let us finish on a more positive note. Although  a derivative need not be continuous, it is  {\sl Darboux Baire -1\/}. That is: (a)  it takes any value
between any two assumed values--- the intermediate value property of continuous functions; see \ref{Br, p.5; R, p.164}; (b)  it is the limit of a
sequence of continuous functions since
 $$
f'(x) = \lim_{n\to \infty} \phi_n(x),\quad {\rm  where}\quad \phi_n(x)=n\left(f(x+\inv{n})- f(x)\right),
$$
 and each $ \phi_n$ is continuous since $f $ is, being differentiable.

 Such functions  have lots of points of  continuity, in fact there are whole intervals of continuity on
every perfect set; \ref{B-B-T, pp.22--23}.  So that on every perfect set the points of discontinuity are nowhere-dense.
\sectn{References}

\napa{[B-D-P] B Bongiorno, L Di Piazza \& D  Preiss}  A constructive minimal integral which includes the  Lebesgue integrable functions and
derivatives, {\sl\ J.\ London Math.\ Soc.},\enspace (2)\ 62\ (2000),\ 117--126.
\napa{[B] N Bourbaki} {\it  Fonctions d'une Variable R\'eelle (Th\'eorie \'El\'ementaire), Chap.I--III\/},  Hermann \& Cie., Paris,  1949.
\napa{[B-B-T]\ A Bruckner, J Bruckner \& B Thomson} {\it Real Analysis\/}, Prentice-Hall, New Jersey, 1997.
\napa{[Br]\ A Bruckner} {\it Differentiation of Real Functions,\/} CRM Monograph Series,Vol.5,
American Mathematical Society, 1994.
\napa{[Bu1]\ P S Bullen} An unconvincing counterexample,{\sl\ Int.\ J.\  Math.\ Educ.\ Sci.\
Technol.},\enspace 19\ (1988),\ 455--459.
\napa{[Bu2]\ P S Bullen} Integration and trigonometric series, to appear.
\napa{[Bus]\ E Busko}   Une relation entre la d\'erivabilit\'e \`a droite et la
continuit\'e, {\sl\ Ens.\ Math.},\enspace 12 (1960)\ 243--247.
\napa{[D]\ A Denjoy} {\it M\'emoire sur la D\'erivation et son Calcul Inverse\/},\ Gauthier-Villars,\ 
Paris, 1954.
\napa{[D-S]\ J D DePree \& C W Swartz }{\it  Introduction to Real Analysis\/}, John Wiley \& Sons, New York, 1988.
\napa{[E]} {\it Encyclopedia of Mathematics, 1--10, Suppl.I\/}, Kluwer Academic Publishers, \enspace Dordrecht, 1988--1997
\napa{[G]\ Alan D. Gluchoff}Trigonometric series and theories of integration, {\sl\
Math.Mag.},\enspace\hfil 67\ (1994),\ 3--20.
\napa{[Go]\ R A Gordon}  {\it The Integrals of Lebesgue, Denjoy, Perron and Henstock\/},  Amer.\
Math.\ Soc Memoir, 1994. 
\napa{[Gr]\ J V Grabiner} {\it The Origins of Cauchy's Rigorous Calculus\/},  MIT Press, Cambridge, 1981
\napa {[H] T Hawkins}  {\it Lebesgue's Theory of Integration. Its Origins and Development\/}, University of Wisconsin Press,
Madison,1970.
\napa{[He]\  R Henstock} {\it Theory of Integration\/}, {\sl\ Butterworth's  London,}\ 1963.
\napa{[Ho]\ E W Hobson} {\it The Theory of Functions of a Real Variable 
and the Theory of Fourier's Series I, II\/}, Cambridge University Press, 1926.
\napa{[J]\ R L Jeffery} {\it The  Theory of Functions of a Real Variable\/}, University of Toronto Press, Toronto, 1953.
\napa{[K]\ H Kestelman} {\it Modern Theories of Integration\/}, Dover Publ., New York.1960.
\napa{[K-K]\ R Kannan \& C K Krueger}{\it  Advanced Analysis on the Real Line\/},  Springer-Verlag, New York, 1996.
\napa{[Kn]\ K Knopp} {\it Theory and Application of Infinite Series\/}, Blackie \& Son Ltd., \ London,\
1948. 
\napa {[Ku]\ J Kurzweil} {\it Nichtabsolut Konvergente Integrale\/},   Leipzig, 1988.
\napa{[L]\ H  Lebesgue} {\it  Le\c cons sur l'Int\'egration et la Recherche des
Fonctions Primitives\/},  Gauthiers--Villars,  Paris; 1st Ed. 1904; 2nd Ed.1928.
\napa{[L-V] Lee Peng   Yee \&  R V\'yborn\'y} {\sl The Integral: An Easy Approach after Kurzweil and Henstock\/}, Aust.  Math.\ Soc. Lecture Series
\#14, Cambridge University Press,  Cambridge, 2000.
\napa{[Mcl]\ R M McLeod} {\it The Generalized Riemann Integral\/}, Carus Math\ Monograph \#
20, Math.\ Assoc,\ America,1980.
\napa{[McS]\ E J McShane} {\it Unified Integration\/}, Academic Press Inc., New York, 1983.
\napa{[N] I P Natanson} {\it Theory of Functions of  Real Variable, I, II,\/}\noter{{\sevenrm Engl. transl by L F Boron of} {\sevencyr I P Natanson,
Teoria Funktsi{\sevenrm\u{\sevencyr i}} Veshchestvenno{\sevenrm\u{\sevencyr i}} Peremenno{\sevenrm\u{\sevencyr i}}}} Frederick Ungar
Publishing Co., New York, 1964
\napa{[P] I N Pesin} {\it Razvite Ponyatiya Integrala\/}\noter{\sevencyr  I N Pesin, Razvite Ponyatiya Integrala.} ,
Moscow,\enspace 1966. {\sl Engl.\
 transl.}\enspace : {\it Classical and Modern Integration Theory\/},\noter{\sevenrm Take care with this translation, some  usages are 
historical,  rather than  modern.} New York, 1970.
\napa{[R] K A Ross} {\it Elementary Analysis: The Theory of Calculus\/}, Springer-Verlag, New York, 1980.
\napa{[Ru] W Rudin} {\it Real and Complex Analysis\/}, McGraw-Hill Book Co., 1966.
\napa{[S]\ S Saks}  {\it Theory of the Integral\/},  2nd Ed.\ rev., Hafner,\ New York,\enspace 1937.
\napa{[W]\ W Walter} {\it Differential- und Integral-Ungleichungen\/}, 	Springer-Verlag, Berlin, 1964.

\bye